\documentclass[preprint,12pt]{amsart}
\usepackage{amssymb}
\usepackage{amsthm}
\usepackage{amsmath}
\usepackage{amsfonts}
\usepackage{stackrel}
\usepackage{appendix}
\usepackage{enumitem}
\usepackage{mathrsfs}
\usepackage{color}
\usepackage[dvipsnames]{xcolor}
\usepackage{babel}
\usepackage[T1]{fontenc}
\usepackage{tikz}
\usepackage{subcaption} 

\newtheorem{thm}{Theorem}[section]

\newtheorem{pro}[thm]{Proposition}

\newtheorem{rem}[thm]{Remark}

\numberwithin{equation}{section}
\usepackage{hyperref}
\hypersetup{colorlinks=true,linkcolor=blue,filecolor=blue, citecolor=blue,urlcolor=cyan}
\usepackage{comment}

\newcommand{\R}{\mathbb{R}}
\newcommand{\N}{\mathbb{N}}

\newcommand{\D}{\displaystyle}

\def\qq#1{\qquad \mbox{#1}\quad}
\def\qt#1{\quad \mbox{#1}\ }
\newcommand{\al}{\alpha}
\newcommand{\be}{\beta}
\newcommand{\De}{\Delta}
\newcommand{\de}{\delta}

\newcommand{\e}{\varepsilon}

\newcommand{\la}{\lambda}
\newcommand{\La}{\Lambda}

\newcommand{\na}{\nabla}
\newcommand{\om}{\omega}
\newcommand{\Om}{\Omega}
\newcommand{\Omb}{\overline{\Om}}
\newcommand{\p}{\partial}

\newcommand{\te}{\theta}

\newcommand{\vf}{\varphi}

\newcommand{\wto}{\rightharpoonup}

\title[Bifurcation  for a quasilinear problem]{Bifurcation for   indefinite weighted $p$-laplacian problems  with  slightly subcritical nonlinearity}
\author{Mabel Cuesta}
\address[M.~Cuesta]{Department of Mathematics, Universit\'e du Littoral C\^ote d'Opale (ULCO), Laboratoire de Math\'ematiques Pures et Appliqu\'ees Joseph Liouville (LMPA),
62100 Calais   (France)}
\email[M.~Cuesta]{Mabel.Cuesta@univ-littoral.fr}
\author{Rosa Pardo}
\address[R.~Pardo]{Departamento de An\'alisis Matem\'atico y Matem\'atica Aplicada,  Universidad Complutense de Madrid, 28040--Madrid, Spain.}
\email[R.~Pardo]{rpardo@ucm.es}
\thanks{The second author is supported by grants  PID2019-103860GB-I00,  MICINN,  Spain, and by UCM, Spain,  Grupo 920894.}

\begin{document}

\maketitle
\begin{abstract}

We study a superlinear elliptic boundary value problem  involving the $p$-laplacian operator, with  changing sign weights.
The problem has positive solutions bifurcating from the trivial solution set at the two principal eigenvalues of the corresponding linear weighted boundary value problem. The two principal eigenvalues are bifurcation points from the trivial solution set to positive solutions.

Drabek's  bifurcation result applies when the nonlinearity   is of power growth.
We extend Drabek's  bifurcation result to {\it slightly subcritical} nonlinearities. Compactness in this setting is a delicate issue obtained via Orlicz spaces. 
\end{abstract}

\keywords
Positive solutions,  subcritical nonlinearity, changing sign weight,  $p$-Laplacian. Orlicz spaces.

\subjclass[2020]{
35B32, 
35J92, 
35B09, 
47J15. 
}

\section{Introduction}
Our aim is to prove the existence of positive solutions to a Dirichlet problem for a class of quasilinear elliptic equations whose nonlinear term has non  power growth, and  involves indefinite nonlinearities. 
More precisely, we consider 
\begin{equation}
\label{eq:ell:pb} 
-\De_p u
=\la V(x) u^{p-1}+m(x)f(u),  \quad \mbox{in } \Om, \qquad u= 0,  \quad \mbox{on } \partial \Om,
\end{equation}
where $p>1$,  $\Om \subset \R ^N $,  with $ N>p$,   is a bounded, connected open set, with  $C^{1,\alpha}$ boundary $\partial \Om$, $\De_p(u)=\text{div }(|\na u|^{p-2}\na u)$ is the $p$-Laplacian operator, $1<p<\infty$, $\la\in\R$ is a real parameter. The weights   $V\in L^\infty (\Om)$  and $m\in C^1(\overline{\Om})$ both changing sign in $\Om$, and $f \in C([0,+\infty)) $ is {\it slightly subcritical} (see \ref{H2}).
Prototype nonlinearities are the following ones
\begin{equation}\label{f:al}
f(s):=\dfrac{|s|^{p^* -2}s}{\big[\ln(e+|s|)\big]^\be},\qquad
f(s):=\dfrac{|s|^{p^* -2}s}{\Big[\ln\big(e+\ln(1+|s|)\big)\Big]^\be},
\end{equation}
where $p^* :=\frac{Np}{N-p}$ is the   critical Sobolev exponent, and $\be>0$ is a fixed exponent, or even having a smaller perturbation of order $O(|s|^{q -2}s)$ at infinity, with $p<q<p^*.$ 

\medskip

There is a huge amount of literature when  $V\equiv 1$, and $f(s)$ is $p$-superlinear and  grows at infinity as, say, $s^{q-1}$ with $q<p^* .$  The case $p=2$ and $m$ changing sign was first studied by \cite{Alama-Tarantello} for a power-like nonlinearity, next by  \cite{BCN}, among  others, and for the case $p\not=2$ see for instance \cite{Bo-Ta, Ka-Ra-Um, Il}. 
\medskip

We focus in widening the solvability of the quasilinear problem \eqref{eq:ell:pb} to nonlinearities $f$ slightly subcritical. The literature is more scarce when the nonlinearity  is not of subcritical power type at infinity. 
We extend bifurcation results for the $p-$laplacian case and for those nonlinearities. Specifically,  we will assume   the following hypothesis on $f$:
\begin{enumerate}[label=\textbf{(f\arabic*)$_\infty$}]
\item\label{H2}
$f$ is  slightly subcritical: ${\D\lim_{s\to +\infty}}\, \frac{f(s)}{|s|^{p^*-1}}=0$;

\item\label{H4} 
$f$ increasing at infinity and there exists two constant $s_0 >0$ and $c_0>1$ such that
\begin{equation}
\label{derf:f}
\frac{sf(s)}{F(s)}\ge c_0,\qquad \forall s>s_0 ,
\end{equation}
where $F$ is the primitive of $f$ satisfying $F(0)=0$.
\end{enumerate}

\smallskip

\medskip

Using local bifurcation  techniques, we prove the existence of positive solutions to problem \eqref{eq:ell:pb}. Some hypothesis on the behavior of $f$ close to 0 are needed. We will assume  the following classical condition:

\smallskip

\newcounter{num}
\setcounter{num}{3}
\begin{enumerate}[label=\textbf{(f\thenum)$_0$}\addtocounter{num}{1}]
\item
\label{H5}
$f(0)=0,$ there exists $\de_0>0$ such that $f(s)>0$ for $s\in(0,\de_0)$, and $f$ is $p$-sublinear at zero : $\ \displaystyle{\lim_{s\to 0}\frac{f(s)}{|s|^{p-1}}=0}$.
\end{enumerate}

To characterize a principal eigenvalue of a changing sign weight $V\in L^\infty (\Om)$, let us consider the following  eigenvalue problem: 
\begin{equation}\label{egV}
\left\{
\begin{array}{rll}
- \Delta_p \vf &=\lambda V(x)|\vf|^{p-2}\vf & \hbox{ in } \Om,\\
\vf&=0 &\hbox{ on } \partial\Omega .
\end{array}
\right.
\end{equation}
It is known that \eqref{egV} possesses exactly two principal eigenvalues,  denoted by $\lambda_1 (V)$ and $\lambda_{-1} (V)$, characterized by 
\begin{equation}\label{egv}
\lambda_1 (V)= \inf_{u\in S^+}\int_\Om |\nabla u|^p dx, \quad  \lambda_{-1}(V)=-\inf_{u\in S^-}\int_\Om |\nabla u|^p dx,
\end{equation}
where 
$$S^\pm=\left\{u\in W_0^{1,p}(\Om)\,: \int_\Om V(x)|u|^p dx =\pm 1\right\}. $$
Those  are the only eigenvalues  associated to a non-negative eigenfunction and they are simple and isolated (see \cite{cuesta1} and references therein).  

\medskip

Let us denote by  $\vf_{1}(V)$ $\big(\vf_{-1}(V)\big)$ the positive eigenfunction associated to $\la_1 (V)$ $\big(\la_{-1} (V)\big)$ of $L^\infty$-norm equal to 1, and  introduce the following hypothesis on $f$  near 0 and on the weight $m$   :

\setcounter{num}{4}
\begin{enumerate}[label=\textbf{ (f\thenum)$_0$}\addtocounter{num}{1}]
\item
\label{H6}
there exists a constant $C_1>0$ and a continuous function 
$g_0:[0,+\infty)\to[0,+\infty)$  such that for all $\tau>0$, 
\begin{equation}
\label{F:mu}
\left|\frac{f(\tau s)}{f(s)}\right|\le C_1 (1+\tau^{p^*-1})\qq{for all} 0<|s|< 1, 
\end{equation}
and
$
{\D\lim_{s\to 0}}\, \dfrac{f(\tau s)}{f(s)}= g_0(\tau)\,,\quad\text{uniformly for }\tau\text{ on compact intervals}.
$

\noindent Moreover, $g(0)=0$, $g_0(\tau)>0$ for all $\tau>0$, and   
\begin{equation}
\label{a:vf1:<0}
\int_{\Om} m(x)g_0(\vf_{\pm 1} (V))\,\vf_{\pm 1} (V)\, dx  < 0. 
\footnote{Note that the continuity of the function $g_0$ and the fact that  $g_0(st)=g_0(s)g_0(t)$ for all $s, t$ positive,  imply that  $g_0(s)=s^{q-1}$ for some $q\leq p^*$.}
\end{equation}
\end{enumerate}

\bigskip

We prove the existence of two continuum of positive solutions bifurca\-ting from the two principal eigenvalues of \eqref{egV}.
Drabek \cite[Theorem 14.18, p. 189]{Drabek} and Tak\'a\v c--Girg \cite[Proposition 3.5]{Girg-Takac} use Browder–Petryshyn topological degree (\cite{Bro-Pe, Skr}), an extension of the Leray-Schauder degree for monotone mappings. Del Pino--Manasevich \cite[Theorem 1.1]{DelPino-Manasevich} use a homotopy to the case $p=2$.  
Roughly speaking,  bifurcation theorems for the $p$-laplacian case apply for subcritical power type nonlinearities, but do not treat slightly  subcritical nonlinearities such as \eqref{f:al}.

\medskip

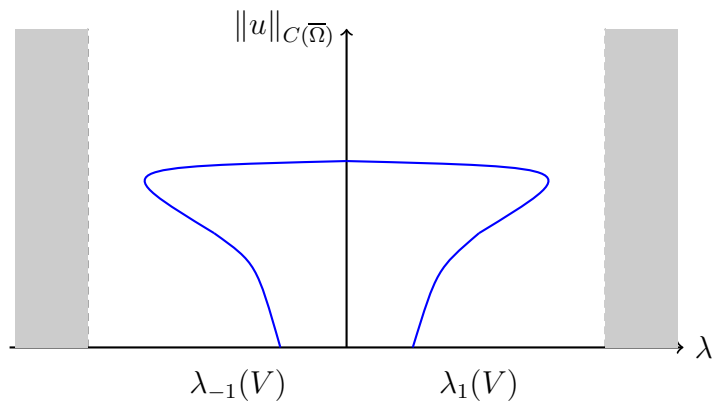
\begin{figure}[h]
\setlength{\unitlength}{.8cm}
\centering 
\begin{tikzpicture}[scale=.8]
\draw[thick,->] (-5.6,0) -- (5.6,0) node[anchor=west] {$\lambda$};
\draw[thick,->] (0,0) -- (0,5.3) node[anchor=east] {$\|u\|_{C(\overline{\Omega})}$}; 
\draw[blue, thick] (1.1,0) .. controls (1.5,1.3) .. (2.2, 1.9);
\draw[blue, thick] (0, 3.1) .. controls (4,3) .. (2.2, 1.9);
\draw[blue, thick] (0, 3.1) .. controls (-4,3) .. (-2.2, 1.9);
\draw[blue, thick] (-1.1,0) .. controls (-1.5,1.4) .. (-2.2, 1.9);
\node[label=below:{$\la_1(V)$}] at (2.2, 0) {};
\node[label=below:{$\la_{-1}(V)$}] at (-1.8, 0) {};
\draw[densely dashed] (4.3,0) -- (4.3, 5.3);
\fill[gray!40!white] (4.3, 0) rectangle (5.5,5.3);
\draw[densely dashed] (-4.3,0) -- (-4.3, 5.3);
\fill[gray!40!white] (-4.3, 0) rectangle (-5.5,5.3);
\end{tikzpicture}
\caption{\small Possible global bifurcation diagram representing a {\it mushroom--shaped} compact component.}
\label{fig}
\end{figure}

When $p=2$, for a subcritical power type nonlinearity, in \cite{Cano-Casanova} the author works on indefinite weights  so that \eqref{eq:ell:pb}  has a bounded, {\it mushroom--shaped} compact continuum of positive solutions connecting  the  two principal eigenvalues of \eqref{egV}, see Fig.\ref{fig}. A  uniform $L^\infty(\Om)$ bound for the solutions is guaranteed through the blow up method. Although we prove some  $L^\infty$  estimates for the solutions to \eqref{eq:ell:pb}, a uniform $L^\infty(\Om)$ bound  is an open problem. \\

Let us  define
\begin{equation*} 
\Om^\pm :=\{x\in \Om : \ \pm m(x)>0\},\qquad 
\Om^0 :=\{x\in \Om : \ m(x)  = 0\},
\end{equation*}
and assume that:
\begin{enumerate}[label=\textbf{(m\arabic*)}]
\item 
\label{m0}
$\sup m^+>0$ and $\sup m^->0$, so both $\Om^\pm \ne\emptyset $;
\item
\label{m1}  let $\om^{+,0}:={\rm int}\,\big(\Om^+\cup \Om^0\big)$, $\om^{+,0}$ have a finite number of connected components, and
$\big|\big\{x\in \om^{+,0}:\, \pm V(x)>0\big\}\big|\not=0$;
\item
\label{m2} let $\om^0: ={\rm int}\,\big(\Om^0\big)$, $\om^0 \not=\emptyset$ have a finite number of connected components, and   $\big|\big\{x\in \om^0:\,\pm  V(x)>0\big\}\big|\not=0$. 
\end{enumerate} 

\smallskip

The pair $(\la, 0)$ is a solution of \eqref{eq:ell:pb}  for every $\la\in\R$, whenever $f(0)=0$. Pairs of this form will be designated as the trivial solutions of \eqref{eq:ell:pb}. We say that $(\la, 0)$ is a {\it bifurcation point} of \eqref{eq:ell:pb}  if in any neighborhood of $(\la, 0)$ there exists a nontrivial solution of \eqref{eq:ell:pb}.
It is well known that if $(\la, 0)$ is a  bifurcation point of \eqref{eq:ell:pb}  then $\la$ is an eigenvalue of \eqref{egV}, see for instance \cite{DelPino-Manasevich}.

Next we state a bifurcation theorem for problem \eqref{eq:ell:pb}.

\begin{thm}\label{th:ex}
Assume that $f$  satisfies hypothesis {\rm {\rm \ref{H2}}}{\textbf -}{\rm \ref{H4}} and also
{\rm \ref{H5}}{\textbf -}{\rm \ref{H6}}.
Assume further that $m\in C^1(\Omb)$ changes sign in $\Om$ and satisfies hypothesis {\rm \ref{m0}}.
Then 
\begin{enumerate}
\item[\rm (i)] there exists two  closed and connected sets $ \mathscr{C}_1^+$ and   $ \mathscr{C}_{-1}^+ $ of  positive  solutions to \eqref{eq:ell:pb},  bifurcating from the trivial solution set at the bifurcation point $(\la_1 (V),0)$ and  $(\la_{-1} (V),0)$ respectively; 

\item[\rm (ii)] either  ${\mathscr C}_1^+$   is unbounded, or 
$(\la_{-1} (V),0)\in {\mathscr C}_1^+.$ Likewise, either   ${\mathscr C}_{-1}^+$   is unbounded, or
$(\la_{1} (V),0)\in {\mathscr C}_{-1}^+$;

\item[\rm (iii)]
$ \mathscr{C}_1^+$ bifurcates to the right and   $ \mathscr{C}_{-1}^+ $ to the left. Precisely,
there exists $\delta>0$ such that 
\begin{equation*}
(\la,u)\in {\mathscr C}_1^+, 0<|\la-\la_1(V)|<\delta, 0\not=\|u\|<\delta \Longrightarrow \la>\la_1(V),
\end{equation*}
and also
\begin{equation*}
(\la,u)\in {\mathscr C}_1^-, 0<|\la-\la_{-1}(V)|<\delta, 0\not=\|u\|<\delta  \Longrightarrow \la<\la_{-1}(V);
\end{equation*}

\item[\rm (iv)] if we assume either {\rm \ref{m1}} or {\rm \ref{m2}},  there exist $\La_1, \La_{-1}\in \R$, with  $$\La_{-1}<\la_{-1}(V)<0<\la_1 (V)< \La_1,$$ such that  if problem \eqref{eq:ell:pb}$_\lambda$ has  non-negative nontrivial solutions, then    $\La_{-1}\le \la \le \La_1$. 
\end{enumerate}
\end{thm}

We will give in Proposition \ref{pro:non:exist} an upper (resp. lower) bound of $\La_1$ (resp. of $\La_{-1}$) .

\bigskip

The main novelties of this paper are the following:

\medskip

(i) We extend bifurcation theorems for $p$-laplacian problem to  slightly  subcritical nonlinearities, see Theorem \ref{pro:bif drabek}. 
Compact embedding  of Sobolev   spaces into Lebesgue spaces (Rellich--Kondrachov Theorem \cite[Theorem 9.16]{Brezis}) do not apply for slightly subcritical functions with non power-like growth. Instead, the compact\-ness  is achieved via Orlicz spaces  under hypothesis {\rm \ref{H2}}, see Proposition \ref{pro:compG} and Theorem \ref{th:comp}.   

(ii) We state sufficient conditions for the existence of a continuum of  positive solutions in a quasilinear bifurcation problem, see Theorem \ref{th:la>la1}, which complements the  results in \cite{DelPino-Manasevich, Girg-Takac}. The branch of positive solutions satisfies Rabinowitz's alternative.

(iii) In Theorem \ref{th:ex} we  generalize partially  results of \cite{Cano-Casanova} on changing sign weights, for $p\ne 2$ and non power-like nonli\-nearities. 

(iv) In Theorem \ref{th:L:infty} we prove some explicit $L^\infty$  estimates for quasilinear problems with a slightly  subcritical non linearity.

\bigskip

The paper is organized in the following way.  In section \ref{sec:non:exist} we prove some explicit $L^\infty$  estimates for quasilinear problems with a slightly  subcritical non linearity.
A Drabek's type bifurcation result    is developed in Section \ref{sec:bif} followed by a a result on  existence of two branches of positive solutions  bifurcating from $\la_1 (V)$ and  $\lambda_{-1}(V)$ respectively, see Theorem \ref{th:la>la1}. In Section \ref{sec:proof}, we  obtain a non 
existence result. We end this Section with the proof of Theorem \ref{th:ex}.
We've tried to  make a summary on tools of Orlicz spaces as simple as possible and self-contained, see Appendix \ref{sec:orlicz}.

\section{$L^{\infty}$ bounds and  classical $C^1$- regularity}

An important issue for  weak  solutions $u$ (positive or not)  of our problem is if there are uniform \textit{ a priori} bounds.
When $\la=0$, $a(x)\equiv 1$ and $f$ is the first case of \eqref{f:al},  the relationship between $L^\infty$, Sobolev and Lebesgue norms  has been studied in \cite{CP2015, CP-ring-like, Castro-Mavinga-Pardo, Cuesta-Pardo_AJM_2022, Mavinga-P-JMAA-17, Pardo_2023} for the semi-linear case, and in  \cite{Damascelli-Pardo, Mavinga-P-JMAA-17, Mavinga-P-MJM-21} for  the $p-$Laplacian case. 
A uniform $L^\infty$ \textit{a priori} bound for  positive weak solution,  is known when $\be> \frac{p}{N-p}$. When $p=2$, it is also known that there is a positive solution blowing up 
as $\be\to 0,$ see \cite{Clapp-Pardo-Pistoia-Saldana} for details.
From a recent result of \cite{Pardo_2024},  first  we  derive an estimate of the $L^\infty$ norm of weak solutions to   problem \eqref{eq:ell:pb} in terms of the $\|\cdot\|_{p^*}$, see Theorem \ref{th:L:infty}.

\smallskip

Throughout this paper, for any $u\in W_0^{1,p}(\Om)$ we will denote 
\label{sec:non:exist}
$$\| u\|:=\left(\int_\Om |\nabla u|^p dx\right)^{1/p}.$$

\subsection{$L^{\infty}$ bounds}

Let us consider  the  quasilinear problem 
\begin{equation}\label{pde}
-\Delta_p u=f(x,u)  \hbox{ in } \Om,\qquad
u=0 \hbox{ on } \partial\Omega,
\end{equation}
and assume that $\Omega$ is a bounded domain of $\R^N$ and   $f:\Omega\times \R\to\R$ is a Carath\'eodory function.

We recall that  $u\in W^{1,p}(\Om)
$ is a   \textit{ weak solution} to \eqref{pde} if  for each  $\psi\in W_0^{1,p}(\Om)$,
\[
\int_{\Om}|\na u|^{p-2}\na u\cdot \na \psi\, dx=  \int_{\Om}f(x,u) \psi \, dx .
\]

If $p<N$ and there exist $0\leq q\leq p^*$ and  $c>0$ \  such that 
\begin{equation}\label{G}
|f(x,s)|\leq c(1+|s|^{q-1}) \qquad \text{a.e. } x\in \Omega,  \quad \text{and for all }  s\in\R, 
\end{equation}
then  any $u\in W_0^{1,p}(\Om)$ solving weakly problem \eqref{pde}, belongs to $L^{\infty}(\Om)$, cf. e.g. \cite{Anane2} in the subcritical case (i.e. $q<p^*$), and \cite{GuVe}  in the critical case (i.e. $q=p^*$).

The following result  is a consequence of \cite[Theorem 1.6]{Pardo_2024}, and gives an estimate of  the $L^\infty$ norm of a solution in terms of its $L^{p^*}$ norm, for nonlinearities $f$  that are  slightly  subcritical in the sense of {\rm \ref{H2}}.

\begin{thm}[$L^\infty $ estimates for slightly subcritical problems]$ $
\label{th:L:infty}
Assume that there exists a continuous function  $\tilde{f}:[0,+\infty)\to [0,+\infty)$ satisfying
\begin{equation}
\label{f:tilde}
|f(x,s)|\leq \tilde{f}(|s|),\quad \forall s\in \R;\quad \tilde{f}(s)>0,\quad \forall  s>0.
\end{equation}
where 
$\tilde{f}$ satisfies {\rm \ref{H2}}.
Let $h(s):=\frac{s^{p^* -1}}{\tilde{f}(s)},$ for $s>0.$
Then, for any $\e >0$ there exists $C=C_\e >0$ such that 
the following holds: 
\begin{equation}\label{H:inf}
h\big(\|u\|_{\infty}\big)\le 
C\ \bigg[\Big(1+\|u\|_{p^*}^{\, \frac{1}{p-1}}\Big)\|u\|_{p^*}\bigg]^{(p^*-1)\frac{p}{N}+\e},
\end{equation}
for every nontrivial weak solution of \eqref{pde}, 
where $C=C(\e,\tilde{f},p,N,|\Om|),$ and it is independent of $u.$ 
\end{thm}

\begin{rem}\label{rem} 1. Note that we can always redefine $\tilde{f}$ in order to be strictly increasing for $s>0$, for instance $\max_{0\leq t\leq s}  \tilde{f}(t) + s^q$, with $p-1<q<p^*-1$.

2. This result also holds for changing sign solutions.    
\end{rem}

\begin{proof}[Proof of Theorem \ref{th:L:infty}]
The proof is adapted from the proof of \cite[Theorem 1.6]{Pardo_2024}, where  it is assumed that
$$
|f(x,s)|\leq |a(x)|\tilde{f}(s),\ \text{a.e }x\in\Om,\ \forall s\in \R;$$
$$ a\in L^r(\Om)\text{ with }N/p<r\le\infty,\ \tilde{f}(s)>0,\ \forall  s\not=0.
$$
Let $u$ be a solution to \eqref{pde}. On the one hand, we  notice that whenever $\tilde{f}$ can be chosen  strictly increasing for $s>0$, then we can always follow the second possibility (ii) in \cite[Theorem 1.6]{Pardo_2024}. 
On the other hand, we point out that in p. 13-14 the following estimate can be read
\begin{equation*}
h\big(\|u\|_{\infty}\big)\le C 
\|a\|_r^{\ \te} \ \Big(\|\tilde{f}(u) \|_{\frac{p^*}{p_{N/r}^*-1}}\Big)^{\te-1}\		
\ \|u\|_{p^*}^{\ \vartheta},
\end{equation*}
for $u=u_k$, where $\{u_k\}$ is a sequence of solutions to \eqref{pde},
$\te =\frac{p_{N/r}^*-1}{p_{N/q}^*-1},$ $\vartheta =\te\ (p_{N/q}^*-p),$ $
\forall q\in\big(N/p,N\big),
$
and
$
p_{N/s}^*:=\frac{p^*}{s'}=p^*\left(1-\frac{1}{s}\right),$ for $1\le s\le\infty,
$
where $s'$ is the conjugate exponent of $s$,  $1/s+1/s'=1$, and $C=C(p,N,q,|\Om|)$.

Hypothesis \eqref{f:tilde} implies that $r=\infty$, and $p_{N/r}^*=p^*$, so
$$
h\big(\|u\|_{\infty}\big)\le C\ \Big(\|\tilde{f}(u) \|_{\frac{p^*}{p^*-1}}\Big)^{\te-1}\		
\ \|u\|_{p^*}^{\ \vartheta},
$$
with 
\begin{align}
\label{def:al:be:inf}
&	\te=\te(q) =\frac{p^*-1}{p_{N/q}^*-1},\qquad
\vartheta=\vartheta(q) =\te\ (p_{N/q}^*-p).
\end{align}
Moreover, from subcriticality (see \ref{H2}),
\begin{equation*}
\|\tilde{f}(u) \|_{\frac{p^*}{p^*-1}}
\le C
\left(1+\|u\|_{{p^*}}^{p^*-1}\right).
\end{equation*}
Consequently
\begin{equation*}
h\big(\|u\|_{\infty}\big)\le C 
\left(1+\|u\|_{{p^*}}^{p^*-1}\right)^{\te-1}\		
\ \|u\|_{p^*}^{\ \vartheta}.
\end{equation*}
Let
$
\Theta(q,t):=\left(1+t^{p^*-1}\right)^{\te-1}\	t^{\, \vartheta}\,.
$
Fixed $t=t_0>0$, the function $q\to \Theta(q,t_0)$ for $q\in \big(N/p,N\big)$,
satisfies the following 
$$
\ln\big(\Theta(q,t_0)\big):=(\te-1)\ln \left(1+t_0^{p^*-1}\right)\	+ \vartheta \ln t_0.
$$
Derivating with respect to $q$, and since $\frac{d}{d\,q}\left(p_{N/q}^*\right)=\frac{p^*}{q^2}$,  
we can write
\begin{align*}
&\frac{\frac{d}{d\,q}\big(\Theta(q,t_0)\big)}{\Theta(q,t_0)}
=\te'(q)\ln \left(1+t_0^{p^*-1}\right)	+ \vartheta'(q)\ln t_0\\
&\quad=-\frac{(p^*-1)}{\big(p_{N/q}^*-1\big)^2}\frac{p^*}{q^2}
\left[\ln \left(1+t_0^{p^*-1}\right)
+\left[\big(p_{N/q}^*-p\big)-\big(p_{N/q}^*-1\big)
\right] \ln (t_0)\right]
\\ 
&\quad=-\frac{(p^*-1)}{\big(p_{N/q}^*-1\big)^2}\frac{p^*}{q^2}
\left[\ln \left(1+t_0^{p^*-1}\right)
-\big(p-1\big) \ln (t_0)\right].
\end{align*}
Consequently,
\begin{equation*}
\frac{d}{d\,q}\big(\Theta(q,t_0)\big)< 0 \iff 
\frac{1}{p-1}> \frac{\ln (t_0)}{\ln \left(1+t_0^{p^*-1}\right)}\,.
\end{equation*}
It is easy to see that the function $t\to\frac{\ln (t)}{\ln \left(1+t\right)}$ is an increasing function for $t>0$, and so
$$
\sup_{t\in\R^+}\frac{\ln (t)}{\ln \left(1+t^{p^*-1}\right)}=\lim_{t\to \infty}\frac{\ln (t)}{\ln \left(1+t^{p^*-1}\right)}
=\frac{1}{p^*-1}< \frac{1}{p-1},
$$
then $\Theta(q,t_0)$ is decreasing in $q$ for any $t=t_0\in\R^+$, so, for any $t$ fixed 
$$
\inf_{q\in (\frac{N}{p},N)} \Theta(q,t)=\Theta\big(N,t\big)=\Big(1+t^{p^*-1}\Big)^{\frac{p'}{N}}\,t^{\frac{p}{N}(p^*-1)}
\le C\bigg[\Big(1+t\Big)^{\frac{1}{p-1}}\,t\bigg]^{(p^*-1)\frac{p}{N}}
\,,
$$
since \eqref{def:al:be:inf}, and then
\begin{align*}
\te (N)-1
&= \frac{p^* -1}{p^* (1-1/N) -1} -1=\frac{p^* /N}{p^* (1-1/N) -1}=\frac{1/N}{1-1/p}=\frac{p'}{N} \\  
\vartheta (N)&=\left(1+\frac{p'}{N}\right)p^*\left(1-\frac{1}{N}-1+\frac{p}{N}\right)=p^*\left(\frac{p-1}{N}+\frac{p}{N^2}\right)\\  
&=\frac{p}{N}p^*\left(1-\frac{1}{p}+\frac{1}{N}\right)=\frac{p}{N}p^*\left(1-\frac{1}{p^*}\right)
=\frac{p}{N}(p^*-1).
\end{align*}
Finally, and since the infimum is not attained in $\big(N/p,N\big)$, 
for any $\e>0$,
there exists a constant $C=C_\e>0$ such that
\begin{equation*}
h\big(\|u\|_{\infty}\big)\le 
C\ \bigg[\Big(1+\|u\|_{p^*}^{\, \frac{1}{p-1}}\Big)\|u\|_{p^*}\bigg]^{(p^*-1)\frac{p}{N}+\e},
\end{equation*}
where  $C$ is independent of $u$, ending the proof.
\end{proof}

\subsection{Regularity of weak solutions} 

Since weak  solutions to problem \eqref{pde}  under condition \eqref{G} belong to $L^\infty (\Om)$, the classical results of \cite{DiBenedetto,  Lieberman, To} for a regular domain $\Omega$ of class $C^{1,\alpha}$ provide that $u$ belongs to
$C^{1,\mu}(\overline{\Om})$ for some $\mu=\mu(N,p,\alpha) \in (0,1)$. Precisely, we have:

\begin{pro}[$C^{1,\mu}(\overline{\Om})$-regularity for Problem \eqref{pde}]
\label{reg1}
Assume  hypo\-thesis \eqref{G}.  Then, for any weak solution  $u\in W_0^{1,p}(\Omega)$  of equation \eqref{pde},  
\begin{enumerate}
\item[\rm (i)] $u\in L^{\infty}(\Om)$ 
and there exists a constant $c_0>0$ depending only  on $p,N, c $ and $\|u\|_{p^*}$ such that
\[
\|u\|_\infty \leq c_0.
\]
\item[\rm (ii)]  There exists $\mu=\mu(p,N,\alpha)\in(0,1)$ for which  $u\in C^{1,\mu}(\overline{\Om})$.
Moreover,   there exists  $C=C(p, N, \alpha, c,\|u\|_{\infty})$ such that 
\[
\|u\|_{C^{1,\mu}(\overline{\Om})}\leq C.
\]
\end{enumerate}
In particular, a  $\, W_0^{1,p}(\Om)$- bounded sequence (or  $L^{p^*}(\Om)$-bounded) of solutions to problem \eqref{eq:ell:pb}$_{\la_n}$, with $\la_n$ varying in a bounded interval, is uniformly bounded in $C^{1,\mu}(\overline{\Om})$.

\end{pro}
\begin{proof} Part (i) and (ii) are well known, see \cite{DiBenedetto, To,  Lieberman}.
We prove the last assertion. Thanks to Theorem \ref{th:L:infty} and Proposition \ref{reg1}, a sequence of solutions to problem \eqref{eq:ell:pb}  $\, W_0^{1,p}(\Om)$- bounded sequence or  $L^{p^*}(\Om)$-bounded,  remains also bounded in $C^{1,\mu}(\overline{\Om})$.
Indeed, since  Theorem \ref{th:L:infty},  $h\big(\|u_n\|_{\infty}\big)\le C$, and due to $h(s)\to\infty$ as $s\to \infty$, hence $\|u_n\|_{\infty}\le C$. Finally, by part (ii), $\|u_n\|_{C^{1,\mu}(\overline{\Om})}$  remains also bounded.   
\end{proof}
\begin{rem}
\label{rem:reg}
In the critical case $q=p^*$, the $\|\cdot\|_{\infty}$-norm    of weak solutions can not be bounded  in terms of the  the Sobolev norm, as was observed by \cite[page 725]{DeF-Go-Ub}. Consequently, a sequence of solutions to problem \eqref{pde} with $q=p^*$, bounded in $W_0^{1,p}(\Om)$, not necessarily remains bounded in $C^{1,\mu}(\overline{\Om})$. 
\end{rem}

\section{A bifurcation result for a quasilinear problem with a slightly subcritical nonlinearity}  \label{sec:bif}

The first main result in this section, Theorem \ref{pro:bif drabek},  proves a slight generali\-zation of  a Drabek's bifurcation theorem for quasilinear equations, see \cite[Theorem 14.18]{Drabek}. These bifurcation results are based on the well known  works  of  \cite[Theorem 1.3]{R71} and \cite[Lemma 1 and Theorem 2]{Dancer}  for the semilinear case. The second main result, Theorem \ref{th:la>la1},  proves Rabinowitz's alternative for a branch of positive solutions.
\medskip

Consider the quasilinear elliptic problem with parameter $\lambda$
\begin{equation}\label{P:la}
\tag{$P_\lambda$}\quad  \left\{
\begin{array}{rll}
-\De_p u &=\ \la  V(x)|u|^{p-2}u+g(\la ,x,u), & \mbox{in } \Om,\\
u&=\ 0,  &\mbox{on } \partial \Om,
\end{array}
\right.    
\end{equation}
where $g:\R\times\Omega\times  \R\to\R$ is
a Carath\'eodory function such that, for any bounded set $J\subset \R$ there exists some continuous  function  $\bar{f}:\R^+\to\R^+$ such that 
\begin{equation}\label{crois g}
|g(\la ,x, s)|\leq 
\bar{f}(|s|)\quad a.e. \, x\in \Omega, \quad \forall (\lambda,s)\in J\times \R,
\end{equation}
where $\bar{f}$ is slightly subcritical. We denote  by $\bar{F}$ the primitive of $\bar{f}$ vanishing at $0$.
\medskip

Let us define the following operators $N,S: W_0^{1,p}(\Omega)\to W^{-1,p'}(\Omega)$ and  $G :\R\times W_0^{1,p}(\Omega)\to W^{-1,p'}(\Omega)$ :
\begin{align}
\label{NSG}
\langle N(u),\psi\rangle &=\int_\Omega |\nabla u|^{p-2}\nabla u\cdot \nabla \psi \, dx,
\\
\langle S(u),\psi\rangle &=\int_\Omega V(x)| u|^{p-2}u \psi \, dx,
\nonumber\\
\langle G(\lambda, u),\psi\rangle &=\int_\Omega g(\la ,x, u)\psi\, dx.\nonumber 
\end{align}
Define the product space $E:=\R\times W_0^{1,p}(\Omega)$ endowed with the  norm $\|(\lambda, u)\|_E:=|\lambda|+\|u\|$.\\

In order to  reformulate Drabek's bifurcation theorem for the equation $N-\lambda S -G(\lambda,\cdot)=0$ we need to assure that $G$ is a {\em compact} operator. 

\begin{pro}\label{pro:compG}
Assume that $g$ satisfies the growth condition \eqref{crois g} with $\bar{f}(0)=0$, 
$\bar{f}$   continuous, strictly increasing  for $s\geq 0$, and satisfies {\rm \ref{H2}}--{\rm \ref{H4}} with $f=\bar{f}$. 

Then for any bounded closed interval $I\subset \R$,  the operator $G:I\times W_0^{1,p}(\Om)\to W^{-1,p'}(\Omega)$ is compact.
\end{pro}

\begin{proof}
Let $\{\la_n\}_{n\in\N}$ be a sequence in $I$, and $\{u_n\}_{n\in\N}$ a bounded sequence of $W_0^{1,p}(\Om)$. Then, up to a subsequence, there exist $\lambda_0\in I$ and $u$ such that $\lambda_n\to\lambda_0$, $u_n \wto u$ weakly in $W_0^{1,p}(\Omega)$, strongly in $L^{p^*-1}$, a.e. and in measure.  Let us show that  the sequence $\{ G(\lambda_n, u_n)\}_{n\in\N}$ converges in $W^{-1,p'}(\Omega)$. 
For any $\psi\in W_0^{1,p}(\Omega)$,
$$
\left|\int_\Omega \big( g(\la_n,x, u_n)-g(\la_0 ,x, u)\big)\psi\, dx\right|\leq \|\psi\|_{p^*}\,\big\|g(\lambda_n, \cdot, u_n)-g(\lambda_0, \cdot, u)\big\|_{(p^*)'}
$$
Let $z_n:=g(\lambda_n, \cdot,u_n)$ and $z:=g(\lambda_0, \cdot,u)$ and let us apply
Theorem \ref{th:comp}  of  Appendix A  to estimate $\|z_n-z\|_{(p^*)'}$.  

We choose  $ a=\bar{f}^{-1}$ and $b(t)=\frac{p^*}{p^* -1} t^{\frac{1}{p^*-1}}$. Since $\bar{f}$ is strictly increasing,  
$a^*(t)=a^{-1}(t)=\bar{f}(t)$,
see  definition \eqref{def:m*}, \eqref{m*:inv}, 
moreover $B(t):=t^{(p^*)'}=t^{\frac{p^*}{p^*-1}}$, for $t\geq 0$.
Remark that, by {\rm \ref{H2}},   $$a^*(s)= \bar{f}(s)\leq cs^{p^*-1} +d, \quad\quad  A^* (s)=\bar{F}(s)\leq \frac{c}{p^*}s^{p^*} +ds$$
for some $c,d>0$.
\medskip
Now, we  check   the  hypotheses of Theorem \ref{th:comp}.\medskip 

1.   The sequence $\{z_n\}_{n\in\N}$ is bounded in $K_A(\Omega)$, i.e.,  for all $n\in \N$, $\displaystyle\int_\Omega A\Big(\big|z_n (x)\big|\Big)\, dx\leq C$ for some $C>0$. Indeed,   using  \eqref{=:M:M*}, 
\begin{align*}
& \int_\Omega A\Big(\big|z_n (x)\big|\Big)\, dx\leq 
\int_\Om A \big(a^*(|u_n (x)|)\big) \, dx \\
&\qquad = \int_\Om  \big( |u_n (x)|a^*(|u_n (x)|)  + A^*(|u_n|)\big)dx
\leq \int_\Om   (c'|u_n|^{p^*} +d')dx\leq C.
\end{align*}
One proves similarly that $z\in K_A (\Omega)$.\medskip

2.  $A$ satisfies the $\De_2$-condition at infinity if there exist constants $k_0>1$ and  $t_0 > 0$ such that, 
\begin{equation*}
\frac{ta(t)}{A(t)} \le k_0, \qq{for} t>t_0,
\end{equation*}
see Proposition \ref{cnys:De2}.

Through the change of variable $s=\bar{f}^{-1}(t)=a(t)$, since \eqref{=:M:M*}  and by   hypothesis {\rm \ref{H4}}, we can write
\begin{align*}
\frac{ta(t)}{A(t)}
=\frac{s\bar{f}(s)}{\bar{f}(s)s -\bar{F}(s)}
\le \frac{c_0 }{c_0-1}=:k_0, \qq{for all} t\ge  \bar{f}(s_0 ).
\end{align*}

\medskip

3.  It remains to prove that $B$ increases essentially more slowly than $A$   at $+\infty$. 
It is equivalent to prove that $A^*$ increases essentially more slowly than $B^*$, 
see \cite[Lemma 13.1]{Krasnoselski-Ruticki}. That is to check, that $\bar{F}$ increases essentially more slowly than $t\to t^{p^*}$, which  is a consequence of 
hypothesis {\rm \ref{H2}}. This concludes the proof.
\end{proof}

\medskip

We define ${\mathscr S}$ as the closure in $E$ of the  non trivial  weak solution set  of \eqref{P:la} :
\begin{equation}\label{defS}
{\mathscr S}:=\overline{\big\{(\lambda, u)\in E :\, u\not=0\,\hbox{ and }  N(u)-\lambda S(u)-G(\lambda, u)=0\big\}}^E.\end{equation}

Now we can prove the following bifurcation theorem :
\begin{thm}\label{pro:bif drabek}
Let $g$ satisfies the growth condition \eqref{crois g} with $\bar{f}(0)=0$, 
$\bar{f}$   continuous,  strictly increasing,  
slightly subcritical at infinity,  see 
{\rm \ref{H2}}, satisfying 
also {\rm \ref{H4}} and  {\rm \ref{H5}} with $f=\bar{f}$. 
Let $\mu_1$ be either $\lambda_1(V) $ or $\lambda_{-1} (V)$.

Then there exist  a maximal closed connected  set ${\mathscr C}$ of ${\mathscr S}$ containing $(\mu_1,0)$  such that  either
\begin{enumerate}
\item[\rm (i)] ${\mathscr C}$ is unbounded on $E$, or else
\item[\rm (ii)] ${\mathscr C}$ contains  a point $(\tilde{\lambda},0)$ where $\tilde{\lambda}\not=\mu_1$  is an eigenvalue of problem \eqref{egV}.
\end{enumerate}
\end{thm}

\begin{proof}
We adapt the proofs of  \cite[Theorem 14.18 and Theorem 14.20]{Drabek},  and of \cite[Proposition 3.5]{Girg-Takac}, based on the well known result of Rabinowitz \cite[Theorem 1.3]{R71} to the case of a nonlinearity $g$ satisfying \eqref{crois g}.
Notice  that   the operator $G$ is compact according to Proposition \ref{pro:compG}. 

Let $deg\,$ denotes the  Browder–Petryshyn topological degree 
(\cite{Bro-Pe, Skr})  defined for  monotone mappings  generalizing the Leray-Schauder  degree, see also  \cite[\S 14.6]{Drabek}.
Since $\mu_1$ is an isolated eigenvalue, reasoning as is \cite[Theorem 14.18]{Drabek}, for any $r>0,$ and 
$\delta>0$ small enough, the following degree is well defined in any open neighborhood $B_r$ of $0$,  and 
$$
deg\, (N-(\mu_1+\delta) S, B_r, 0)\not= deg\, (N-(\mu_1- \delta) S, B_r, 0).
$$ 
Assume for a moment that
\begin{equation}\label{in0}
\lim_{\|u\|\to 0} \frac{\|G(\lambda,u)\|_{W^{-1,p'}(\Omega)}}{\|u\|^{p-1}}=0,
\end{equation}
uniformly for $\lambda$ on bounded sets. Then, the invariance for homotopy 
yields  that  for some $r,\ \delta>0$ small enough,
\begin{align}
\label{deg:ne}
&deg\, (N-(\mu_1+\delta) S -G(\mu_1+\delta,\cdot), B_r, 0)\\
&\qquad \ne deg\, (N-(\mu_1- \delta) S-G(\mu_1-\delta,\cdot), B_r, 0) ,   \nonumber
\end{align}
ending the proof.

To prove \eqref{in0}, let  $0\not=\{u_n\}_{n\in\N}$ be a sequence in $W_0^{1,p}(\Om)$ such that $\|u_n\|\to 0$. Then, for $\lambda$ on a bounded set $J\subset\R$
$$
\lim_{n\to\infty}\frac{\|G(\lambda,u_n)\|_{ W^{-1,p'}(\Omega)}}{\|u_n\|^{p-1}}=
\lim_{n\to\infty }\sup_{ \|\psi\|\leq 1}\int_{\Om} \frac{|g(\la ,x, u_n)\psi| }{\|u_n\|^{p-1}}\, dx
$$
Let us denote $w_n=\frac{u_n}{\|u_n\|}$, fix $\e >0$ and  $\psi\in W_0^{1,p}(\Om)$ with $\|\psi\|=1$. 
Since {\rm \ref{H5}},  there exists $\delta_1 >0$ such that, for $\la\in J$,
$$ 
|g(\la ,x, s)|\leq\e |s|^{p-1}, \quad \text{a.e.}\ x\in \Om,\ \forall |s|\leq \delta_1.
$$
Denote for all $\delta>0$, 
$$
\Om_n^{\delta}:=\{x\in\Om : \, |u_n(x)|\geq \delta\}.
$$
Since  we can assume that $u_n\to 0$ in measure,  for any $\delta >0$ fixed $|\Om_n^{\delta} |\to 0$. So, there exists an $n_0=n_0(\e )\big(=n_0(\e ,\delta_1)\big)$ such that 
$|\Om_n^{\delta_1}|\le \e $ for all $n\ge n_0$.

Here and thereafter $C$ stands for a constant independent of $n$.  By Holder's inequality,
\begin{equation}\label{in1}
\int_{\Om\setminus \Om_n^{\delta_1}} \frac{|g(\la ,x, u_n)\psi| }{\|u_n\|^{p-1}}\, dx\leq  \e \, \| \psi\|_p\,   \|w_n\|_p^{p-1}\leq C \e .
\end{equation}
Using  that $\bar{f}$ satisfies {\rm \ref{H2}}, there exists $\delta_2>\delta_1$ such that, 
$$|g(\la ,x, s)|\leq \e |s|^{p^*-1}, \quad\text{a.e.}\ x\in \Om,\  \forall |s|\geq \delta_2.$$
Besides, let $c(\delta_1,\delta_2)>0$ be such that 
$$
\bar{f}(|s|)\leq c(\delta_1, \delta_2) |s|^{p-1},\quad \forall \delta_1\leq |s|\leq \delta_2
$$
and let $n_1=n_1(\e )\big(=n_1(\e ,\delta_1,\delta_2)\big)\in \N$ be such that 
$$
|\Om_n^{\delta_1}|^{\frac{p^*-p}{p^*}}\leq \frac{\e }{c(\delta_1,\delta_2)}, \quad \forall n\geq n_1.
$$
Then, since 
$W_0^{1,p}(\Om)$ is embedded into $L^{p^*} (\Om)$, using Holder's inequality, we have,
\begin{align}\label{in2}
\int_{\Om_n^{\delta_1}} \frac{|g(\la ,x, u_n)\psi| }{\|u_n\|^{p-1}}\, dx &\leq \e  \int_{\Om_n^{\delta_2}} \frac{|u_n|^{p^* -1} |\psi|}{\|u_n\|^{p-1}}\, dx + c(\delta_1, \delta_2)\int_{\Om_n^{\delta_1} \setminus \Om_n^{{\delta_2}}}  |w_n|^{p-1}\,|\psi|\, dx \nonumber 
\\
&\leq  C\, \e \|u_n\|^{p^* -p}\|\psi\|_{p^*} + c(\delta_1, \delta_2)\|w_n\|_{p^*}^{p-1} \|\psi\|_{p^*} |\Om_n^{\delta_1}|^{\frac{p^*-p}{p^* }}\\
&\leq C\e +C c(\delta_1, \delta_2) |\Om_n^{\delta_1}|^{\frac{p^*-p}{p^* }}
\leq C\e ,\quad  \forall n\geq n_1.\nonumber
\end{align} 
Inequalities \eqref{in1} and \eqref{in2}  give \eqref{in0}.
\end{proof}

We complete this section by a regularity result for solutions bifurcating from zero, that we will use later.

\begin{pro}\label{pro:vf11}
Assume the hypotheses  of Theorem \ref{pro:bif drabek}. 
Let $\{(\lambda_n, u_n)\}_{n\in\N}$ be  a sequence of solutions to \eqref{P:la} for $\la=\la_n$, such that  $\displaystyle\lim_{n\to\infty}\lambda_n=\la^*$. Then the following three
statements are equivalent, as $n\to\infty$:
\begin{itemize}
\item[\rm (i)]\ $\|u_n\|_{p^*}\to 0$;
\item[\rm (ii)]\ $\|u_n\|_{\infty}\to 0$;
\item[\rm (iii)]\ $\|u_n\|_{C^{1,\mu}(\overline{\Om})}\to 0$.
\end{itemize}
Moreover, the following hold
\begin{enumerate}
\item[\rm (1)] if $\displaystyle\lim_{n\to\infty}\lambda_n=\la_{\pm 1} (V) $ , and $\|u_n\|_{p^*}\to 0$, then  the sequence $w_n:=\frac{u_n}{\|u_n\|_{\infty}}$ is the union of two disjoint sub-sequences $\{w_n^{'}\}_{n\in\N}$
and $\{w_n^{''}\}_{n\in\N}$, one
of them possibly empty, such that, if nonempty, they satisfy $w_n^{'}  \to \vf_{\pm 1} (V)$ ,
and/or $w_n^{''} \to -\vf_{\pm 1}(V)$ 
in $C^{1,\mu'}(\overline{\Om}) $  as $n\to\infty$. Here, $\mu' \in (0, \mu)$ is arbitrary and $\mu $ is given in  Proposition \ref{reg1}.
\item[\rm (2)] 
If $\|u_n\|_{p^*}\to 0$,
then
$\displaystyle\lim_{n\to\infty}\lambda_n=\la_{\pm 1} (V)$.  
\end{enumerate}
\end{pro}
\begin{proof}
Clearly, (iii) implies (i) and (ii). We will prove that (i)$\implies $ (ii) and (ii)$\Longrightarrow $ (iii).\\

\noindent  (i)$\Longrightarrow $ (ii).
Our hypothesis  on $g$ guarantees that, for all $n\in\N,$ 
\eqref{P:la} for $\la=\la_n$
satisfies  the hypothesis of   Theorem \ref{th:L:infty}. So,  using \eqref{crois g}, we
apply Theorem \ref{th:L:infty}  with  $\tilde{f} (s):=(|\la_{\pm 1} (V) |+\delta)\|V\|_\infty s^{p-1}+\bar{f} (s)$ for $s>0.$
Fix $\e>0$,  by Theorem \ref{th:L:infty} and their conclusion \eqref{H:inf}, $h(\|u_n\|_\infty) \to 0$ as $n\to\infty$ with $h(s):=\frac{s^{p^* -1}}{\tilde{f}(s)}$ for $s>0.$ Since  $h$ is a continuous function, $h(s)>0,$ $\forall  s>0$, and  
$$
\lim_{s\to 0^+}\frac{1}{h(s)}=\lim_{s\to 0^+} (|\lambda_{\pm 1}(V) | +\delta )\|V\|_{\infty}s^{p-p^*} +\frac{\bar{f}(s)}{s^{p^*-1}}= +\infty,
$$ 
then $h(0)=0$ and $\|u_n\|_\infty\to 0$ as $n\to\infty$, and (ii) holds.\\

\noindent  (ii)$\Longrightarrow $ (iii).  The function $w_n:= u_n /\|u_n\|_{\infty}$ satisfies
\begin{equation}\label{res}\int_\Om |\nabla w_n|^{p-2}\nabla w_n\cdot \nabla \psi \, dx =\int_\Om \kappa_n(x)\psi \, dx
\end{equation}
for all $\psi\in W_0^{1,p}(\Om)$,  with 
$$\kappa_n(x):=
\lambda_n  V(x)|w_n|^{p-2}w_n  + \frac{g(\lambda_n, x, u_n)}{\|u_n\|_{\infty}^{p-1}}.
$$
Since $\bar{f}$ satisfies {\rm \ref{H2}} and   {\rm \ref{H5}}, and  $\|w_n\|_{\infty}=1$,   there exists $D>0$ such that, for all $n$ large enough,
$$
|\kappa_n(x)|\leq D\quad \hbox{ a.e. } x\in \Om, 
$$
and from Proposition \ref{reg1},  there exists a constant $C>0$
$$\|w_n\|_{C^{1,\mu}(\overline{\Om})}\leq C=C(p,N, D)$$
and (iii) follows.\medskip

Now we prove the last two parts of the lemma.\\ 

(1) The embedding $C^{1,\mu}(\overline{\Om})\to C^{1,\mu'}(\overline{\Om})$ being compact, by Arzel\`a-Ascoli's theorem, the sequence $\{w_n\}_{n\in\N}$
contains a subsequence that
converges into  $C^{1,\mu'}(\overline{\Om})$ to some $w$; we denote it again by $w_n \to w$. Notice that $\|w\|_\infty = 1$.  We let $n\to \infty$ in equation \eqref{res}, to conclude that $w \in C^{1,\mu'}(\overline{\Om})$ must satisfy
$$\int_\Om |\nabla w|^{p-2}\nabla w\cdot \nabla \psi \, dx =\la_{\pm 1}(V) \int_\Om  V(x)|w|^{p-2}w \psi \, dx.$$
By the simplicity of $\lambda_{\pm 1} (V)$, 
we have $w=\pm \vf_{\pm 1} (V).$\\

(2) Assume  that $\la_n\to \la^*$.
Then, reasoning as in part (1), $0< w_n\to w$ in $C^{1,\mu'}(\Omb)$
and $w\ge 0$, $w\ne 0$, satisfy
$$\int_\Om |\nabla w|^{p-2}\nabla w\cdot \nabla \psi \, dx =\la^* \int_\Om  V(x)|w|^{p-2}w \psi \, dx,$$
so $w$ is a non-negative eigenfunction, and necessarily $\la^*=\la_{\pm 1}(V) $. 
\end{proof}

\begin{thm}\label{th:la>la1}
Let $g$ satisfies the growth condition \eqref{crois g} with $\bar{f}(0)=0$, 
$\bar{f}$   continuous,  strictly increasing,  
slightly subcritical at infinity,  see 
{\rm \ref{H2}}, satisfying 
also {\rm \ref{H4}} and  {\rm \ref{H5}} with $f=\bar{f}$.

Denote $E=\R\times W_0^{1,p}(\Om)$ and 
$$
{\mathscr P}^+:=\left\{u\in C_0^1(\overline{\Om}): u>0 \hbox{ in } \Om \hbox{ and } \frac{\partial u}{\partial \nu}<0 \hbox{ on } \partial \Om \right\}.
$$ 
Then,
\begin{enumerate}
\item[\rm (i)]  there exists a closed connected set $ {\mathscr C}_1^+\subset E$ of  positive weak   solutions to \eqref{P:la}  containing  the bifurcation point $(\la_1 (V),0)$;
\item[\rm (ii)]  there exists a closed connected set $ {\mathscr C}_{-1}^+\subset E$ of  positive weak   solutions to \eqref{P:la}  containing  the bifurcation point $(\la_{-1} (V),0)$;

\item[\rm (iii)]  either  ${\mathscr C}_1^+$   is unbounded, or 
$(\la_{-1} (V),0)\in {\mathscr C}_1^+.$ Likewise, either   ${\mathscr C}_{-1}^+$   is unbounded, or
$(\la_{1} (V),0)\in {\mathscr C}_{-1}^+$;

\end{enumerate}
\end{thm}
\begin{proof}
Let us check that all   the hypothesis  of bifurcation Theorem \ref{pro:bif drabek} are accomplished, for $\hat{g}(\la ,x, s)$, where $\hat{g}(\la ,x, \cdot)$ is the odd extension of $g(\la ,x, \cdot)$.

Consider the modified problem 
\begin{equation}\label{prob-i}
-\Delta_p u =\ \lambda V(x)|u|^{p-2}u +\hat{g} (\la ,x, 	u) \hbox{ in }\Om,\qquad
u =\ 0 \hbox{ on } \partial\Om.
\end{equation}
Trivially, if $u\leq 0$ is a solution of \eqref{prob-i} then $-u$ is a positive solution of \eqref{prob-i},  and positive solutions to  \eqref{prob-i} are positive solutions to \eqref{P:la}.

Denote $\mu_1$ either $\la_1(V)$ or $\la_{-1}(V)$  and let  $G: E\to W^{-1,p}(\Om)$ be the map 
$$
\langle G(\la,u),\psi\rangle= \int_\Om \hat{g}(\la ,x, 	u)\psi\, dx.
$$ 
Theorem \ref{pro:bif drabek}  applied to the   principal  eigenvalue $\mu_1=\la_1(V)$, implies the existence of a maximal closed connected    set  $\hat{{\mathscr C}}_1$  (with respect to  the closure of nontrivial solutions ${\mathscr S}$) for \eqref{prob-i}, containing $(\la_1(V),0)$.
Likewise, whe get  the existence of another maximal closed connected  set $\hat{{\mathscr C}}_{-1}$  containing $(\la_{-1}(V),0)$.  Both  sets $\hat{{\mathscr C}}_{\pm 1}$ satisfy properties (i) and (ii) of Theorem \ref{pro:bif drabek}.  Moreover,  by oddity of $\hat{g}$, if $(\la,u)\in \hat{{\mathscr C}}_1,$ also $(\la,-u)\in \hat{{\mathscr C}}_1$, in other words $-\hat{{\mathscr C}}_{1} =\hat{{\mathscr C}}_{ 1}$  and   $-\hat{{\mathscr C}}_{-1} =\hat{{\mathscr C}}_{-1}$. \\

\medskip
\noindent(i) Let us  first prove that there exists $\delta>0$ such that
\begin{equation}\label{Si:P}
\hat{{\mathscr C}_1}\cap B_\delta \big((\lambda_1 (V) ,0)\big) \setminus \big\{(\lambda_1 (V),0)\big\} \subset (\lambda_1 (V)-\de, \la_1 (V)+\de)\times \left({\mathscr P}^+\cup -{\mathscr P}^+ \right).
\end{equation} 
Assume by contradiction that there exists  a changing sign sequence of solutions $(\lambda_n, u_n)\in \R\times( W_0^{1,p}(\Om)\setminus \{0\})$ such that
$$\int_\Om |\nabla u_n|^{p-2}\nabla u_n \cdot \nabla \psi \, dx=\lambda_n \int_\Om  V(x)|u_n|^{p-2}u_n \psi \, dx +\int_\Om \hat{g}(\la ,x, 	u_n) \psi \, dx$$
for all $\psi \in W_0^{1,p}(\Omega) $ and $(\lambda_n, u_n)\to (\la_1 (V),0).$  Denote by $w_n:=\frac{u_n}{\|u_n\|_\infty}$.  Due to  the previous Proposition \ref{pro:vf11},
\begin{equation}\label{vn:to:vf1}
w_n:=\frac{u_n}{\|u_n\|_\infty}\to \pm\vf_1 (V), \qq{in} C^{1,\mu'}(\Omb)\,.
\end{equation}
Since   $\vf_1 (V)\in {\mathscr P}^+$ and ${\mathscr P}^+$  is an open set in the $C_0^1$-topology, there exists a neighbor\-hood $\mathscr{V}\subset {\mathscr P}^+$ such that  $w_n\in \mathscr{V}$, and we conclude that also $u_n\in {\mathscr P}^+$, contradicting that  $\{u_n\}$ is a changing sign sequence of solutions. If $w_n\to -\vf_1 (V)$, we take instead the sequence of solutions    $\{-w_n\}_{n\in\N}$ and reason as in  the previous case.  We have proved \eqref{Si:P}.

Notice that $\hat{{\mathscr C}_{1}}\cap\big(\R\times {\mathscr P}^+\big)\ne\emptyset$  as a result of \eqref{Si:P} and the fact that $\hat{{\mathscr C}}_1=-\hat{{\mathscr C}_1}$. 
Let $\hat{{\mathscr C}_1}^+$ be the maximal subcontinuum of ${\mathscr S}$ containing the branch of positive  
solutions bifurcating from $(\la_{1} (V),0)$  contained in $\hat{{\mathscr C}_{1}}$. 
By a subcontinuum of ${\mathscr S}$ we mean a subset of ${\mathscr S}$  which is closed and connected in the Banach space $E$. 

Taking ${\mathscr C}_{1}^+:=\hat{{\mathscr C}_{1}}^+\ne\emptyset$, let us prove that 
\begin{equation}
\label{posi}
{\mathscr C}_{1}^+ \subset \R\times {\mathscr P}^+\cup \{(\la_{\pm 1}(V),0)\}.\end{equation}
Assume by contradiction that there exits $(\la^*,u^*)\in {\mathscr C}_{1}^+$ with $(\la^*,u^*)\not=(\la_{\pm 1}(V),0)$ and that  there exits  $\{(\la_n, u_n)\}_n$ in ${\mathscr C}_{1}^+ \cap \R\times {\mathscr P}^+$ such that 
$\la_n\to\la^*, u_n\to u^*$ en $W_0^{1,p}(\Om)$. Since $u_n> 0$ then $u^*\geq 0$. If $u^*=0$ then $(\la^*,0)$ is a bifurcation point of positive solution and therefore $\la^*=\la_{\pm 1}(V)$, a contradiction. Since  $u^*\geq 0, u^*\not=0$ then $u^*$ is a nontrivial  non-negative  weak solution to problem \eqref{P:la} for $\la=\la^*$, and therefore, by the regularity results of Proposition \ref{reg1} and the standard strong maximum principles  for quasilinear elliptic equations  (\cite{Pucci-Serrin, Pucci-Serrin-Zou, Vazquez}) we deduce that $u^*>0$ in $\Om$ and  
$u^*\in C^{1,\mu}(\overline{\Om})$, so the claim \eqref{posi} is proved.\\

\noindent (ii) The proof is similar to the case (i) with minor changes.\\ 

\noindent (iii) 
To prove that either ${\mathscr C}_1^+$      is unbounded, or $(\la_{-1} (V),0)\in {\mathscr C}_1^+$,
we use $\hat{{\mathscr C}_{1}}^+$. We claim that 
\begin{equation}\label{nc}
\hat{{\mathscr C}}_1\subset  \hat{{\mathscr C}_{1}}^+ \cup \big(-\hat{{\mathscr C}_{1}}^+ \big)\cup \{(\la_{-1}(V),0)\}.
\end{equation}
Since \eqref{posi}, the branch $\hat{{\mathscr C}_{1}}^+\setminus \{(\la_{\pm 1}(V),0)\}$ only contains strictly positive solutions.

In particular,  a pair $(\la^*, u^*)\in{\mathscr C}_1^+$ can not be approached by changing sign solutions, and  the claim \eqref{nc} is proved.
Consequently,   recalling that  $\hat{{\mathscr C}}_1$  satisfies Rabinowitz's alternative, either  ${\mathscr C}_1^+=\hat{{\mathscr C}_{1}}^+ $  is unbounded or else $(\la_{-1}(V),0)\in {\mathscr C}_1^+$.
\end{proof}

\section{Non Existence and Proof of Theorem \ref{th:ex} }
\label{sec:proof}

We will made  use of the well known {\bf Picone's identity} \cite{Allegretto-Huang}, that we recall here. Let us  denote
\begin{align}
\label{L}
&L(v_1,v_2):= |\na v_1 |^p +(p-1) \frac {v_1^p} {v_2^p} |\na v_2 |^p - 
p\frac {v_1^{p-1}}{v_2^{p-1}} \na v_1 \cdot  |\na v_2 |^{p-2} \na v_2, \\
\nonumber 
&R(v_1,v_2):= |\nabla v_1|^p -|\nabla v_2|^{p-2}\nabla v_2 \cdot \nabla\left( \frac{v_1^p}{v_2^{p-1}}\right),\end{align}
for any couple of $a.e.$ differentiable functions $v_1$ and $v_2$ defined on an open set $\mathcal{O}\subset\R^N$, with $v_2>0$.
Then  one has :\\
1. $0\leq L(v_1,v_2)=R(v_1,v_2)$.\\
2. If $v_1/ v_2 \in W_{loc}^{1,1}(\mathcal{O})$ then $L(v_1,v_2)=0$  a.e. in ${\mathcal O}$  implies that  $v_1=cv_2$ for some constant $c\in \R$ in each connected component of $\mathcal{O}$.

\subsection{A  non-existence result}
Let $\om\subset\Om$ be an open bounded  set, possibly non-connected, with $\{\om_i\}_{i\in I}$ its family of open connected components.  We recall that if $\Omega'\subsetneq\Omega$ in an open  domain, then
$ \la_{\pm 1}(V,\Omega')>\la_{\pm 1} (V,\Omega)$, see \cite{cuesta1}. Let us denote by 
$$\la_{1}(V,\om)=\inf_{i\in I}\la_1 (V,\om_i), \quad  
\la_{-1}(V,\om):=\sup_{i\in I}\la_{-1}(V,\om_i) $$
where $\la_{\pm 1} (V,\om_i)$ are resp. the  smallest positive eigenvalue and largest negative eigenvalue of the Dirichlet eigenvalue problem 
\begin{equation*}\label{eq:eigen}
-\De_p \vf=\la V(x)\vf^{p-1} \qt{in}\om_i,\qquad \vf= 0\qt{on}\p \om_i, 
\end{equation*}
see \eqref{egv}. Clearly, if $I$ is finite,  then $\la_{\pm 1} (V,\om)=\la_{\pm 1} (V,\om_i)$ for some $i$  (not necessarily unique). Observe that if 
\begin{equation*}
\label{condV+}
\big|\{x\in\om_i : V(x)  > 0 \}\big| >0 ,\qq{then} \la_1\big(V,\om_i\big)<+\infty.
\end{equation*}
Likewise, if
\begin{equation*}
\label{condV-}
\big|\{x\in\om_i : V(x)  < 0 \}\big|>0,\qq{then} \la_{-1}\big(V,\om_i\big)>-\infty.
\end{equation*}
Let us   agree that $\la_{1}(V,\om_i)=+\infty $  if $V\leq 0$ a.e. in $\om_i $ (resp.   $\la_{-1}(V,\om_i)=-\infty $ if $ V\geq 0$ a.e. in $\om_i$).
Finally, when $V\equiv 1$ we will    write $\la_1(\om)=\la_1(1,\om)$, and 
$\la_{\pm 1}(V)=\la_{\pm 1}(V,\Om)$, when $\om=\Om$.
\begin{pro}
\label{pro:non:exist}
Let  $f$ satisfy hypothesis   {\rm \ref{H2}} and {\rm \ref{H5}}. 
\begin{enumerate}
\item[\rm (i)] Assume that $m$ satisfies {\rm \ref{m0}}  and  {\rm \ref{m1}}. Set $\om^{+,0}:={\rm int}\,\big(\Om^+\cup \Om^0\big)$.
Let us denote $\alpha_{+,0}:= 1+\frac{C_0\sup m^+}{\la_1 ( \om^{+,0})},$ and
\begin{align*}
C_0&:=\inf\{C\geq 0:\, f(s)+Cs^{p-1} > 0 \hbox{ for all } s> 0\}<+\infty .
\end{align*}
If problem   \eqref{eq:ell:pb}$_\la$ has a non-negative nontrivial solution, then
$$
-\infty<\alpha_{+,0}\,\la_{-1}\big(V,\om^{+,0}\big) < \la<\alpha_{+,0}\,\la_{1}\big(V, \om^{+,0}\big)<+\infty.$$
\item[\rm (ii)] Assume  hypothesis {\rm \ref{m2}}.
Set $\om^0: ={\rm int}\,\big(\Om^0\big)$. If problem   \eqref{eq:ell:pb}$_\la$ has a non-negative nontrivial solution, then
$$
-\infty<\la_{-1}\big(V,\om^{0}\big)<\la<\la_1\big(V,\om^{0}\big)<+\infty.
$$
\end{enumerate}
\end{pro}

\begin{proof}

\noindent (i) 
Let $\la>\alpha_{+,0}\la_1\big(V, \om^{+,0}\big)$ and assume by contradiction that there exists a non-negative nontrivial solution $u\in  W_0^{1,p} (\Om)$ to   \eqref{eq:ell:pb}$_\la$. By  standard maximum principles for quasilinear elliptic equations (see \cite{Pucci-Serrin, Pucci-Serrin-Zou, Vazquez}) we deduce that $u>0$ in $\Om$ and $\frac{\partial u}{\partial \nu}<0$ on $\partial\Om$.

Let us write
$\om^{+,0}=\cup_{i}^{m}w_i$ the decomposition of $\om^{+,0} $ as the disjoint union of its open connected components.
Then $\la_{1}(V, \Om^{0,+})=\la_1 (V,\om_i)$  for some connected component $\om_i$ of $\om^{+,0}$.  
Let $\vf_1^{+,0}:=\vf_1\big(V,\om_i\big) $ be a positive eigenfunction associated to $\la_1(V,\om^{+,0})$ of $L^\infty$-norm equal to $1$.   
For simplicity, we will also denote by $\vf_1^{+,0}$ the  extension by $0$ of $\vf_1^{+,0}$ in all $\Om$. 

Fix $\e >0$ and put $\xi:=\frac{(\vf_1^{+,0})^p}{(u+\e)^{p-1}}$. 
If we  multiply   the equation \eqref{eq:ell:pb}$_\la$  by $\xi$ and  integrate  along $\Om$ we find, 
\begin{align*}
&\int_{\Om}|\na u|^{p-2}\na u\cdot \na \left(\frac{(\vf_1^{+,0})^p}{(u+\e)^{p-1}}\right)\, dx
=  \la\int_{\Om}V(x)(\vf_1^{+,0})^{\, p}\frac{u^{p-1}}{(u+\e)^{p-1}}\\
\nonumber 
&\qquad\qquad\qquad\qquad +\int_{\Om} m^+(x)f(u)\frac{(\vf_1^{+,0})^p}{(u+\e)^{p-1}} \, dx .
\end{align*}

Using Picone's identity  we deduce

\begin{align*}
0	&\le\int_{\Om}L\big(\vf_1^{+,0}, u+\e \big)\, dx=\int_{\Om}R\big(\vf_1^{+,0}, u+\e \big)\, dx\\
\nonumber 
&= \int_{\Om}\big|\na (\vf_1^{+,0})\big|^{p} \, dx -\int_{\Om}|\na u|^{p-2}\na u\cdot \na \left(\frac{(\vf_1^{+,0})^p}{(u+\e)^{p-1}}\right)\, dx\\
\nonumber
& = \int_{\Om}\big|\na (\vf_1^{+,0})\big|^{p} \, dx -\la\int_{\Om}V(x)(\vf_1^{+,0})^{\, p}\left(\frac{u}{u+\e}\right)^{p-1}\\
\nonumber
& \qquad -\int_{\om^{+,0}} m^+(x)f(u)\frac{(\vf_1^{+,0})^p}{(u+\e)^{p-1}} \, dx\\
\nonumber
&\leq \int_{\om^{+,0}}\big|\na (\vf_1^{+,0})\big|^{p} \, dx -\la\int_{\om^{+,0}}V(x)(\vf_1^{+,0})^{\, p}\left(\frac{u}{u+\e}\right)^{p-1}\\
\nonumber
& \qquad 
+C_0\sup \, m^+\int_{\om^{+,0}}\left(
\frac{u}{u+\e}\right)^{p-1}(\vf_1^{+,0})^{p} \, dx .
\end{align*}

Letting $\e\to 0$, we obtain 
$$0\leq \lim_{\e\to 0}\int_{\Om} L(\vf_1^{+,0},u+\e)\, dx\le\left(\alpha_{0,+} -\frac{\la}{\la_1 (V,\om^{+,0})}\right)\int_{\om^{+,0}}|\na (\vf_1^{+,0})|^{p} \, dx$$
which is clearly  a contradiction whenever $\la>\alpha_{+,0}\la_1(V,\om^{+,0})$.
\medskip
In the case $\la=\alpha_{+,0}\la_1(V,\om^{+,0})$ we have 
$$\lim_{\e\to 0} \int_{\Om} L(\vf_{ 1}^{+,0}, u+\e)\, dx  =0
$$
and  using Fatou's lemma
$$ \int_{\Om} \liminf_{\e \to 0} L(\vf_{1}^{+,0}, u+\e)\, dx  \leq 0$$
Since $L(\vf_{ 1}^{+,0}, u+\e)\ge 0$ for all $\e>0$ then we get,
$$ 
L(\vf_{1}^{+,0}, u)= \liminf_{\e \to 0} L(\vf_{1}^{+,0}, u+\e) =0 \quad a.e. 
\, \hbox{ in }\, \Om.
$$
Since trivially 
$\frac{\vf_1^{+,0}}{u}\in W_{loc}^{1,1}(\Om)$, from  the results of Picone's identity 
mentioned above we infer that, in $\om_i$,   $\vf_{ 1}^{+,0}=c u$    for some $c\in\R$. But $\vf_{1}^{+,0}$ vanishes on $\p\big(\om_i\big)\cap\Om $
although  $u>0$ in $\Om$,   a contradiction.

\medskip

(ii)   The proof of non existence  for $\la\leq \la_{-1}\big(V,\om^{+,0}\big)$ is  similar  to the previous case and we omit it.
\end{proof}

\subsection{Proof of Theorem \ref{th:ex}.}
\label{subsec:proof}
Since hypotheses on $f$,  choosing  $\bar{f}$ strictly increasing as in Remark \ref{rem}, hypotheses of Theorem \ref{th:la>la1} are satisfied. Hence, parts (i)--(ii) are accomplished.

\bigskip

\noindent (iii) We have to prove that there exists a neighborhood $\mathscr{U}\subset \R\times W_0^{1,p}(\Om)$ of $(\la_1 (V),0)$ such that  for all $(\la,u_\la)$  in ${\mathscr C}^+\cap \mathscr{U}$, $\la>\la_1 (V)$.
Assume that there is a sequence $(\la_n,u_n)$ of bifurcated positive solutions to  \eqref{eq:ell:pb} in a neighborhood $\mathscr{U}\subset \R\times W_0^{1,p}(\Om)$ of $(\la_1 (V),0)$,   with $(\la_n,u_n)\to(\la_1 (V),0)$.
Let $\e>0$ be fixed. 
Using Picone's identity, see definition \eqref{L},  and that $(\la_1(V),\vf_1(V))$ is an eigenpair, we deduce that for all $\e>0$
\begin{align}
\label{pc}
0	&\leq \int_\Om L( u_n, \vf_1(V)+\e)\, dx =\int_\Om R (u_n, \vf_1 (V)+\e)\, dx \\
\nonumber 
&\leq  \int_{\Om}|\na u_n|^{p} \, dx -\int_{\Om}|\na \vf_1 (V)|^{p-2}\na \vf_1 (V)\cdot \na \left(\frac{u_n^p}{(\vf_1 (V)+\e)^{p-1}}\right)\, dx\\
\nonumber
& = \int_{\Om} m(x)f(u_n) u_n\, dx + \int_{\Om} V(x)u_n^{\, p} \bigg[\la_n-\la_1 (V) \left(\frac{\vf_1 (V)}{\vf_1 (V)+\e}\right)^{p-1}\bigg]\, dx
.
\end{align}
Since  hypothesis {\rm \ref{H6}} on $f$ and \eqref{vn:to:vf1}, we deduce for $w_n=\frac{u_n}{\|u_n\|_\infty}$ that 
$$
\D\frac{f\big(w_n(x)\|u_n\|_\infty\big)}{f\big(\|u_n\|_\infty\big)}\to g_0(\vf_1 (V))\qq{a.e. as} n\to \infty. 
$$
Fix $\de\in(0,1)$.   Hypothesis {\rm \ref{H6}} on $f$, implies that for all $\e_1>0$ there exists a $n_1>0$ such that 
\begin{align*}
\left|\frac{f\big(\mu\|u_n\|_\infty\big)}{f\big(\|u_n\|_\infty\big)}- g_0(\mu)\right|\le \e_1,\ \forall n>n_1, \ 
\forall \mu\in [(1-\de)\vf_1 (V),(1+\de)\vf_1 (V)].
\end{align*}
Clearly $(1-\de)\vf_1 (V)\le w_n\le (1+\de)\vf_1 (V)$ for $n$ big enough, and so
\begin{align*}
\frac{f\big(w_n(x)\|u_n\|_\infty\big)}{f\big(\|u_n\|_\infty\big)}- g_0(w_n(x))\to 0\qq{as}n\to\infty.
\end{align*}
Hence,  using \eqref{vn:to:vf1}, we deduce
\begin{align*}
\int_{\Om}m(x)\frac{f\big(u_n(x)\big)}{f\big(\|u_n\|_\infty\big)}\frac{u_n }{\|u_n\|_\infty}\, dx 
&\underset{ n\to \infty}{\to} \int_{\Om}m(x)g_0(\vf_1 (V))\,\vf_1 (V)\, dx=:I_1\,<0\,,
\end{align*} 
by hypothesis {\rm \ref{H6}}.

Now we prove that $\la_n>\la_1 (V).$ Indeed, on the one hand
$$
\lim_{n\to\infty}\frac{1}{\|u_n\|_p^p}\int_\Om V(x)u_n^p\, dx =\int_\Om V(x)\vf_{1}(V)^p \, dx =\frac{\|\vf_{1}(V)\|^p}{\la_{1}(V)}>0
$$
and then $\int_\Om V(x)u_n^p dx >0$ for $n$ large enough. On the other hand, letting $\e \to 0$ in \eqref{pc}, dividing  by $\|u_n\|_\infty f(\|u_n\|_\infty)$ and using  \eqref{vn:to:vf1}
\begin{align*}
&\frac{\big(\la_1 (V)-\la_n\big)}{\|u_n\|_\infty f\big(\|u_n\|_\infty\big)}\int_{\Om} V(x) u_n^p\le \int_{\Om}m(x) \frac{f(u_n)}{f\big(\|u_n\|_\infty\big)}\frac{u_n }{\|u_n\|_\infty}\, dx\underset{ n\to \infty}{\to} I_1<0 \ .
\end{align*}
We have proved the first part of (iii).  The second part is identical to that one.\\

\noindent (iv) Let us finally define 
\begin{equation*}
\La_1:=\sup\big\{\la>0 : \,  \eqref{eq:ell:pb}_\la \ \text{admits a positive solution} \big\},
\end{equation*}
\begin{equation*}
\La_{-1}:=\inf\big\{\la<0 : \,  \eqref{eq:ell:pb}_\la \ \text{admits a positive solution} \big\}.
\end{equation*}
Since the above, $\La_{-1}<\la_{-1}(V)<\la_1 (V)<\La_1$. Moreover, as a consequence of Proposition \ref{pro:non:exist},  $\La_1$ ,  $\La_{-1}$ are finite.
We have proved  (iv). 
\qed

\begin{appendices}

\section{The convex set  $K_A (\Omega)$ and its gauge $\|\cdot\|_A$ }
\label{sec:orlicz}

\renewcommand{\thesection}{\Alph{section}}
\numberwithin{equation}{section}

Let us summarize some  definitions and results on convex sets and its Minkowski functional (gauge). We refer the reader  to  \cite{Adams, Dacorogna, Krasnoselski-Ruticki}.\medskip

\noindent\textit{ 1.  The convex set $K_A(\Om)$}\\
Let $a :\R^+\to \R^+ $ be  increasing   (that is, $0\leq s\leq t$ implies $a(s)\leq a(t)$) and   righ continuous. 
Let $ A :\R^+\to \R^+$  be defined as
\begin{equation}\label{def:M}
A(t)=\int_0^t 
a(s)\,ds.
\end{equation}
Then $ A$  is a continuous,  strictly increasing and  convex function such that  $A(0)=0$. 

If furthermore $a(0)=0$ and $\D\lim_{s\to+\infty}a(s)=+\infty$, we will say that \textbf{ $A$ is a $N$-function}.

\medskip

Let $\Om$ be an open subset of $\R^N$. The {\textbf convex set  $K_A(\Om)$} is defined as the set of real-valued measurable functions $u:\Om\to\R$ such that :
$$
\int_\Om A(|u(x)|) \, dx< +\infty .\quad 
$$
For all $u\in K_A(\Omega)$ one defines the {\it gauge}  
\begin{equation}\label{def:norM}
\|u\|_{A}:=\inf\left\{\la>0\, : \int_\Om A\left(\frac{|u(x)|}{\la}\right)\, dx\leq 1 \,\right\}.
\end{equation}
Observe that  $\|s u\|_A=|s|\|u\|_A$  for all $s\in [-1,1]$,
$\|u+v\|_A\leq \|u\|_A +\|v\|_A$ for all $u,v\in K_A (\Omega)$  and $\|u\|_A=0$ if and only if $u=0$ a.e.
Note also that $K_A(\Om)$ is a convex set which  is not a linear space in general.\footnote{ The map $u\to \|u\|_A $ is  the Minkowski functional of the convex set
$C =\{u\in K_A(\Omega):\, \int_\Om A(|u(x)|) \, dx\leq  1\}$.} From  the convexity of $A$  it follows that 
\begin{equation}\label{+1}
\forall u\in K_A(\Om), \quad 	\|u\|_{A}\leq  \max\left\{\int_\Om A(|u(x)|)\, dx, \, 1\right\}.\quad
\end{equation}
Furthermore one have 
\begin{equation}\label{pr}
\forall u\in K_A(\Om), u\not=0, \quad 
\int_\Omega A\left(\dfrac{|u(x)|}{\|u\|_A}\right)\, dx \leq 1.
\end{equation}

\noindent\textit{ 2. Conjugate of a $N$- function, generalization of Young's and Holder's inequalities}

Let $A$ be a $N$-function. The  {\textbf conjugate $A^*$ of  $A$} is defined by
\begin{equation}\label{def:M*}
A^* (t):=\max\{st-A(s),\,  s>0\}, \quad \forall t\geq 0
\end{equation}
or, equivalently, 
\begin{equation}\label{def:M*2}
A^* (t):=\int_0^t a^*(s)\, ds, \quad \forall t\geq 0, 
\end{equation}
where $a^*:\R^+\to\R^+$ is given by 
\begin{equation}\label{def:m*}
a^*(t):=\sup\{ s \, : a(s)\leq t\}.
\end{equation}
Observe that for any $t\geq 0$, $a^*\big(a(t)\big)\ge t$ and,  whenever $a$ is  strictly increasing, 
\begin{equation}\label{m*:inv}
a^* (t)=  a^{-1} (t). 
\end{equation}
If follows directly from definition \eqref{def:M*2} that $A^*$ is a $N$-function as well.
It follows also directly  from definition \eqref{def:M*}, the following {\textbf Young's inequality}  :
\begin{equation}\label{def:M:M*}
\forall s,t\in \R^+, \quad st\leq A(s)+ A^*(t).
\end{equation}
Equality holds in \eqref{def:M:M*} if and only if $t=a(s)$ or $s=a^* (t)$. In particular $\forall s,t\in \R^+, $
\begin{equation}\label{=:M:M*}
sa(s)= A(s)+ A^*(a(s)),\quad ta^* (t)= A(a^* (t))+ A^*(t).
\end{equation}

\smallskip

Let $\Omega\subset \R^N$ be  open. Using \eqref{def:M:M*} and \eqref{pr} the following {\bf generalization of Holder's inequality} holds : 
\begin{equation}\label{HI}
\forall u\in  K_A (\Om),\, \forall v\in K_{A^*}(\Om), \quad \left|\int_\Om uv\, dx\right|\leq 2 \|u\|_{A}\, \|v\|_{A^*}.
\end{equation}

\noindent\textit{ 3. $\Delta_2$-condition and comparison of  convex functions}

Let $A$ be a $N$-function, $A$ is said to satisfy the  {\bf $\Delta_2$-condition} if for all $r>1$ there exists $k=k(r)>0$ such that :
\begin{equation}\label{delta2}
A(rt)\leq k A(t) \quad \forall t\ge 0.
\end{equation}
When equation \eqref{delta2} holds only for all $t\ge t_0$, for  some $t_0\ge 0$, then $A$ is said to  satisfy the $\Delta_2$-condition {\it  near infinity}.

\begin{pro}
\label{cnys:De2}
Let $A$ be a $N$-function. A necessary and sufficient condition for $A$  to satisfy the $\De_2$-condition at infinity is 
that there exist constants $k_0>1$ and  $t_0 > 0$ such that, 
\begin{equation*}
\frac{ta(t)}{A(t)} \le k_0, \qq{for} t>t_0.
\end{equation*}    
\end{pro}
See \cite[Theorem 4.1]{Krasnoselski-Ruticki}. 
\medskip

Let $A$ and $B$ be two $N$-functions. We say that  \textbf{ $B$ increases  essentially more  slowly  } than  $A$ \textit{ near infinity}  if 
\begin{equation}\label{<<} 
\forall \delta > 1\qquad  \displaystyle\lim_{t\to+\infty} \frac{B(\delta t)}{A(t)}=0.\end{equation}

We have the following  result:
\begin{thm}\label{th:comp}
Let $\Om\subset \R^N$ be an open set with finite volume and  let $A$ and $B$ be two N-functions. Assume also that 
\begin{enumerate}
\item[\rm (1)] $A$ satisfies the $\Delta_2$-condition near infinity,
\item[\rm (2)] $B$ increases essentially  more slowly than $A$ near infinity. 
\end{enumerate}
Let $\{u_n\}_{n\in \N}$ be a $\|\cdot\|_{A}$-bounded sequence   of $K_A (\Omega)$ converging in measure to some $u\in K_A(\Omega)$. Then    
$$
\lim_{n\to\infty}\|u_n-u\|_B=0.
$$
\end{thm}
Starting from boundedness and convergence in measure, this theorem is a simplified approach to a certain strong convergence. Alternatively, one could use compact embedding in Orlicz-Sobolev spaces \cite{Cianchi}.
\begin{proof}[Proof of Theorem \ref{th:comp}]
Let $0<\e <2\sup_j\|u_j\|_A$ be fixed  and put $v_n:=\frac{|u_n-u|}{\e } $. 
Let us show that there exists $K=K(\e )>0$ such that 
\begin{equation}\label{vn:K}
\forall n\in \N,\quad  \|v_n\|_A \leq K .   
\end{equation}
Indeed, let us choose $r:=\max\Big\{\frac{2 \sup_j \|u_j\|_A}{\e },\frac{2}{\e }\Big\} $,  and  $k=k(r)$, $t_0>0$ as  given in the $\Delta_2$-condition, then 
\begin{equation}\label{ee}
\forall t\geq 0,\quad A(rt)\leq A(rt_0) + kA(t)\end{equation}
Since $A$ is increasing, convex  and satisfies\eqref{ee}, we have
\begin{align*}
A(v_n)&\leq \frac{1}{2}A\left(\frac{2\|u_n\|_A}{\e }\frac{|u_n|}{\|u_n\|_A}\right)+\frac{1}{2}A\left(\frac{2|u|}{\e }\right)    \\
&\leq \frac{1}{2}\left( 2A(rt_0)+ k A\left(\frac{|u_n|}{\|u_n\|_A}\right)+ kA(|u|) \right).
\end{align*}
Hence, using \eqref{+1} and \eqref{pr}, 
\begin{align*}
\|v_n\|_A &\leq \int_\Omega A(|v_n(x)|)\, dx +1 \\
&\leq  \frac{1}{2}\left(2A(rt_0)|\Omega|+ k +k\int_\Omega A (|u(x)|)\, dx \right)+1 :=K,
\end{align*}
and so \eqref{vn:K} is sastisfied. Since \eqref{<<}, let us choose $t_1>0$   such that
$$ B(t)\leq \frac{1}{4}A(t/K)\quad \forall t\geq t_1.$$
Set   
$$
\Omega_n:=\left\{x\in \Omega\, : v_n(x)>B^{-1}\left(\frac{1}{2|\Omega|}\right)\right\}, \
\Omega_n^{'} =\{x\in \Omega_n \, : v_n(x)\geq t_1\}.
$$
Since the sequence $\{u_n\}_n$ converges in measure to $u$, and $\e$ is fixed, $\{v_n\}_n$ converges in measure to $0$ for $\delta=\frac{1}{4B(t_1)}$ there exists $n_0 \in \N$ such that 
$\forall n\geq n_0,\ |\Omega_n|\leq \delta.$
Thus 
\begin{align*}
\int_\Omega\vspace{-.2cm}  B(v_n(x))\, dx 
&=\int_{\Omega\setminus\Omega_n}\vspace{-.2cm} B(v_n(x))\, dx 
+ \int_{\Omega_n^{'}}\vspace{-.2cm}  B(v_n(x))\, dx 
+ \int_{\Omega_n\setminus \Omega_n^{'}}\vspace{-.2cm}  B(v_n(x))\, dx\\
&\leq \frac{|\Omega|}{2|\Omega|} + \frac{1}{4} \int_{\Omega_n^{'}}A\left(\frac{v_n(x)}{K}\right)\, dx +\delta B(t_1)\leq 1,
\end{align*}
thanks to \eqref{vn:K}.
Consequently, $\int_\Omega\vspace{-.2cm}  B\left(\frac{|u_n(x)-u(x)|}{\e}\right)\, dx\le 1 ,$ and so
$$\forall n\geq n_0, \quad \|u_n-u\|_B\leq \e .$$
\end{proof}
\end{appendices}

\section{Acknowledgments}

This work was done during a series of  Pardo's  visits to the Universit\'e du Littoral C\^ote d'Opale ULCO, and of   Cuesta's  visits to the Universi\-dad Complutense de Madrid UCM, whose invitations and hospitality they thanks. The second author is supported by grants PID2019-103860GB-I00, and PID2022-137074NB-I00,  MICINN,  Spain, and by UCM-BSCH, Spain, GR58/08, Grupo 920894.

\end{document}